\documentclass[3p]{elsarticle}

\usepackage{amsmath}
\usepackage{amssymb}
\usepackage{amsthm}
\usepackage[all]{xy}

\newtheorem{lemma}{Lemma}[section]

\newtheorem{thm}{Theorem}[section]

\newcommand{\lag}{\mathfrak g}
\newcommand{\lap}{\mathfrak p}
\newcommand{\laq}{\mathfrak q}

\newcommand{\lasl}{\mathfrak{sl}}
\newcommand{\lagl}{\mathfrak{gl}}

\newcommand{\R}{\mathbb R}

\newcommand{\Sp}{\mathbb S}
\newcommand{\Z}{\mathbb Z}
\newcommand{\C}{\mathbb C}
\newcommand{\V}{\mathcal V}
\newcommand{\CP}{\mathbb{CP}}

\newcommand{\U}{\mathcal U}
\newcommand{\osu}{\mathcal U}
\newcommand{\osv}{\mathcal V}

\newcommand{\osw}{\mathcal W}
\newcommand{\pt}{\mathcal P}
\newcommand{\osm}{\mathcal M}
\newcommand{\vbc}{\mathcal C_\lambda}
\newcommand{\sO}{\mathcal O}

\newcommand{\SO}{\textrm{SO}}
\newcommand{\G}{\textsc G}
\newcommand{\rgr}{\textsc{V}^{\mathbb R}_0(2,n+4)}
\newcommand{\cgr}{\textsc{V}^{\mathbb C}_0}
\newcommand{\LGP}{\textsc P}
\newcommand{\mV}{\textsc V}
\newcommand{\mW}{\textsc W}
\newcommand{\SL}{\textsc{SL}}
\newcommand{\Spin}{\textsc{Spin}}
\newcommand{\GL}{\textsc{GL}}
\newcommand{\LGR}{\textsc R}
\newcommand{\LGQ}{\textsc Q}

\newcommand{\ra}{\rightarrow}

\begin{document}

\title{Penrose transform and monogenic functions.}
\author{Tom\'{a}\v{s} Sala\v{c}}
\ead{salac@karlin.mff.cuni.cz}
\address{Faculty of Mathematics and Physics, Charles University, Sokolovsk\'{a} 83, 18675 Prague, Czech Republic}
\tnotetext[t1]{Research supported by GACR 201/09/H012 and SVV-2011-263317.}

\begin{abstract} Penrose transform gives an isomorphism between the kernel of a $2$-Dirac operator and sheaf cohomology on the twistor space. The point of this paper is to write down this isomorphism explicitly. This gives a new insight into the structure of the kernel of the operator. \end{abstract}


\maketitle

\section{Introduction.}
Let us denote by $\rgr$ the Grassmanian of null $2$-planes in $\R^{n+2,2}$ with a quadratic form of signature $(n+2,2)$. This space is the homogeneous space of a parabolic geometry. There is a sequence
\begin{equation}\label{sequence of operators}
\xymatrix{\Gamma(\mathcal V_1)\ar[r]^{D_1}&\Gamma(\mathcal V_2)\ar[r]^{D_2}&\Gamma(\mathcal V_3)\ar[r]^{D_3}&\Gamma(\mathcal V_4)\ra0}
\end{equation}
of operators living on $\rgr$ which are invariant with respect to the action of the principal group of the geometry, see \cite{TS}. This sequence belongs to the singular character and to the stable range. The first operator $D_1$ is a first order operator and is called the $2$-Dirac operator (in the parabolic setting). Sections which belong to the kernel of $D_1$ are called monogenic sections. The sequence (\ref{sequence of operators}) is interesting from the point of the $2$-Dirac operator (in the Euclidean setting) studied in Clifford analysis, see \cite{CSSS}.

The sequence (\ref{sequence of operators}) is coming from Penrose transform. The Penrose transform is explained in \cite{BE} and \cite{WW}. The Penrose transform for similar sequences in the unstable range was studied in \cite{K}. Similar sequences were obtained by the Penrose transform also in \cite{B}. The Penrose transform lives in the holomorphic category. We loose no information in holomorphic category about the sequence (\ref{sequence of operators}) since $\rgr$ is algebraic variety and we consider real analytic sections which extend uniquely to holomorphic sections over small holomorphic thickenings of affine open subsets of $\rgr$. 

The Penrose transform tells us that the sequence (\ref{sequence of operators}) is exact (after adding the kernel of $D_1$ at the first spot in (\ref{sequence of operators})) on every convex sufficiently small subset of $\rgr$ and is a resolution of the sheaf of monogenic sections. In particular, the Penrose transform gives an isomorphism between the kernel of the operator $D_1$ on affine subset $\mathcal U$ and $H^3(\osw,\sO_\lambda)$. However, general machinery of the Penrose transform does not yield this isomorphism explicitly. In this paper we write this isomorphism explicitly for the case $n=6$. The explicit formula for the Penrose transform is in (\ref{integrating over fibers}). This will give us some interesting insight into the space of monogenic sections. Moreover we decompose the space of monogenic sections with respect to a maximal reductive subgroup of parabolic subgroup, see Theorem \ref{decomposition}.

\subsection{Notation.}\label{notation} 
Let $\cgr(k,n)$ denotes the space of $k$-dimensional subspaces in $\C^n$ that are totally null with respect to some non-degenerate, symmetric, complex bilinear form. The space $\cgr(k,n)$ is smooth algebraic variety and thus also complex manifold. We will denote by $M(k,n,\C)$, resp. $M(k,\C)$, resp. $A(k,\C)$ the space of complex matrices of rank $k\times n$, resp. $k\times k$, resp. the space of skew-symmetric matrices of rank $k\times k$. We denote the span of vectors $v_1,\ldots,v_k$ by $\langle v_1,\ldots,v_k\rangle$. Affine charts on Grassmannian of $k$-planes in $\C^n$ will be denoted by usual notation. We will always work with some preferred basis $\{v_,\ldots,v_n\}$ of $\C^n$. These affine charts will correspond to affine subspaces in the space of matrices $M(n,k,\C)$. If $A=(a_{ji})_{j=1,\ldots,n}^{i=1,\ldots,k}\in M(n,k,\C)$ is such a matrix then we associate to it the $k$-plane spanned by the vectors $u_1,\ldots,u_k$ where $u_i=\sum\limits_j a_{ji}v_j$.

\section{The parabolic geometry.}
Let $h$ be a non-degenerate, symmetric, complex bilinear form on $\C^{10}$. We will work with basis \newline
$\{e_1,e_2,e_3,e_4,e_5,\bar e_3,\bar e_4,\bar e_5,\bar e_1,\bar e_2\}$ such that for all $1\le i,j\le 5:h(e_i,\bar e_j)=\delta_{ij},h(e_i,e_j)=h(\bar e_i,\bar e_j)=0$. Let $\G:=\{g\in End_\C(\C^{10})|\forall u,v\in\C^{10}:h(gu,gv)=h(u,v),\det(g)=1\}$. Let $x_0:=\langle e_1,e_2\rangle$ and let $\LGP:=\{g\in\G|g(x_0)=x_0\}$. The space $\G/\LGP$ is naturally isomorphic to the Grassmannian $\cgr(2,10)$. Let $\G_0$ be the stabilizer of $x^c_0:=\langle e_3,e_4,e_5,\bar e_3,\bar e_4,\bar e_5\rangle$ in $\LGP$. Then $\G_0$ is maximal reductive subgroup of $\LGP$ isomorphic to $\GL(2,\C)\times\SO(6,\C)$. The Lie algebra $\lag$ of $\G$ consists of all matrices of the form
\begin{equation}\label{lie algebra g}
\left(
\begin{array}{cccc}
A&Y_1&Y_2&Y_{12}\\
X_1&B&D&-Y^T_2\\
X_2&C&-B^T&-Y_1^T\\
X_{12}&-X_2^T&-X_1^T&-A^T\\
\end{array}
\right)
\end{equation}
where $A\in M(2,\C),B\in M(3,\C),C,D\in A(3,\C),Y_i,X_i^T\in M(2,3,\C),X_{12},Y_{12}\in A(2,\C)$. There is $\G_0$-invariant gradation 
$$\lag\cong\lag_{-2}\oplus\lag_{-1}\oplus\lag_0\oplus\lag_1\oplus\lag_2$$ 
such that $\lag_0$ is the Lie algebra of $\G_0$ and the Lie algebra $\lap$ of $\LGP$ is $\lag_0\oplus\lag_1\oplus\lag_2$. The subalgebra $\lag_0$ corresponds to the blocks $A,B,C,D$, $\lag_1$ to the blocks $Y_1,Y_2$, $\lag_2$ to the block $Y_{12}$, $\lag_{-2}$ to the block $X_{12}$ and $\lag_{-1}$ to the blocks $X_1,X_2$. Let us denote by $\lag_-:=\lag_{-2}\oplus\lag_{-1}$.

Let $\pi:\G\ra\G/\LGP$ be the canonical projection. Let $\G_-:=\exp(\lag_-)$ and let $\U:=\pi(\G_-)$. The map $\exp$ is biholomorphism between $\lag_-$ and $\G_-$ and the map $\pi$ is biholomorphism between $\G_-$ and $\osu$. 

\section{The double fibration.}
The Penrose transform starts with two fibrations $\tau,\eta$ which fit into diagram
\begin{equation}\label{double fibration diagram}
\xymatrix{&\ar[dl]^\eta\G/\LGQ\ar[dr]^\tau&\\
\G/\LGR&&\G/\LGP.}
\end{equation}
The twistor space $\G/\LGR$ is the connected component of $z_0:=\langle e_1,\ldots,e_5\rangle$ in the Grassmannian $\cgr(5,10)$. We set $\LGR:=\{g\in\G|g(z_0)=z_0\}$. The Lie algebra of $\LGR$ is the subspace of $\lag$ where the matrices $C,X_2,X_{12}$ from (\ref{lie algebra g}) are zero. Let $\bar z_0=\langle\bar e_1,\ldots,e_5\rangle$ and let $\LGR_0:=\{g\in\LGR|g(\bar z_0)=\bar z_0\}$. Then $\LGR_0$ is maximal reductive subgroup of $\LGR$ isomorphic to $\GL(5,\C)$ and its Lie algebra consists of the matrices where $A,B,X_1,Y_1$ are arbitrary and the other matrices are zero.

The correspondence space $\G/\LGQ$ consists of the pairs $(z,x)$ with $x\in\G/\LGP,z\in\G/\LGR$ such that $x\subset z$. We choose $\LGQ:=\LGP\cap\LGR$ so that the Lie algebra $\laq$ of $\LGQ$ sits in the blocks $A,B,D,Y_1,Y_2,Y_{12}$. Let $\LGQ_0:=\G_0\cap\LGR_0$. Then $\LGQ_0$ is maximal reductive subgroup of $\LGQ$ isomorphic to $\GL(2,\C)\times\GL(3,\C)$. The Lie algebra $\laq_0$ of $\LGQ_0$ corresponds to the blocks $A,B$

Let us denote $\osv:=\tau^{-1}(\osu)$ and let $\osw:=\eta(\osv)$. 

\subsection{Projection $\tau$.}\label{subsection tau} 
Let us recall that $x_0=\langle e_1,e_2\rangle,x^c_0=\langle e_3,\ldots,e_5,\bar e_3,\ldots,\bar e_5\rangle,z_0=\langle e_1,\ldots,e_5\rangle$. Let us notice that $x_0^\bot=x_0\oplus x_0^c$. 
The projection $\tau|_\V:\V\ra\U$ sends $(z,x)\mapsto x$. The fibre $\tau^{-1}(x_0)$ is isomorphic to the set of all null 5-planes in $\G/\LGR$ which contains the null 2-plane $x_0$. Thus $(z_0,x_0)\in\tau^{-1}(x_0)$. It is easy to see that any such null 5-plane is determined by an unique null 3-plane in $x^c_0$. For example $z_0=x_0\oplus y_0$ where $y_0:=\langle e_3,e_4,e_5\rangle$. The fibre $\tau^{-1}(x)$ is isomorphic to $\LGP/\LGQ$ and in particular is connected. We deduce that the fibre is biholomorphic to the connected component of $y_0$ in the space of all null $3$-planes in $x_0^c\cong\C^6$. This connected component is biholomorphic to the family of $\alpha$-planes in $\C^6$.

\subsection{The family of $\alpha$-planes in the Grassmanian $\cgr(3,6)$.} 
The restriction of the $\G$-invariant quadratic form on $\C^{10}$ induced by $h$ descends to non-degenerate quadratic form on $x^c_0\cong\C^6$. It will be convenient to make the following identifications. We identify $x_0^c\cong\Lambda^2\C^4$ such that the basis $\{e_3,e_4,e_5,\bar e_3,\bar e_4,\bar e_5\}$ of $x_0^c$ goes to the basis
\begin{equation}\label{null basis of c6}
\{f_0\wedge f_1,f_0\wedge f_2,f_0\wedge f_3,f_2\wedge f_3,-f_1\wedge f_3,f_1\wedge f_2\},
\end{equation}  
of $\Lambda^2\C^4$ where $\{f_0,f_1,f_2,f_3\}$ is the standard basis of $\C^4$. Under this identification, the quadratic form on $\C^6$ goes to the quadratic form $Q$ on $\Lambda^2\C^4$ which is determined by $\alpha\wedge\alpha=Q(\alpha)f_0\wedge f_1\wedge f_2\wedge f_3$ for all $\alpha\in\Lambda^2\C^4$. This identifies $\SL(4,\C)\cong\Spin(6,\C)$. The corresponding isomorphism $\mathfrak{sl}(4,\C)\cong\mathfrak{so}(6,\C)$ is determined by
\begin{eqnarray}\label{isom sl(4,C) and so(6,C)}
\left(
\begin{matrix}
A_1&E_{12}&0&0\\
E_{21}&A_2&E_{23}&0\\
0&E_{32}&A_3&E_{34}\\
0&0&E_{43}&A_4\\
\end{matrix}
\right)\ra
\left(
\begin{matrix}
A_1+A_2&E_{23}&0&0&0&0&\\
E_{32}&A_1+A_3&E_{34}&0&0&E_{12}\\
0&E_{43}&A_1+A_4&0&-E_{12}&0\\
0&0&0&-A_1-A_2&-E_{32}&0\\
0&0&-E_{21}&-E_{23}&-A_1-A_3&-E_{43}\\
0&E_{21}&0&0&-E_{34}&-A_1-A_4\\
\end{matrix}
\right).
\end{eqnarray}
The Grassmannian $\cgr(3,6)$ is the disjoint sum of two families. The first family, called the family of $\alpha$-planes, can be identified with $\CP^3$ by the following mapping. Let $\pi\in\CP^3$. Let $v\in\C^4$ be a representative of $\pi$. Let $\{v,v_1,v_2,v_3\}$ be a basis of $\C^4$. We assign to $\pi$ the $3$-plane $\langle v\wedge v_1,v\wedge v_2,v\wedge v_3\rangle,i=1,2,3$. It is easy to see that the $3$-plane is null and that the map is well defined. The latter family, called the family of $\beta$-planes, can be identified with $\mathbb P(\C^4)^\ast$ by the assignment $[\omega]\in\mathbb P(\C^4)^\ast\mapsto\langle v_1\wedge v_2,v_1\wedge v_3,v_2\wedge v_3\rangle$ where $\{v_1,v_2,v_3\}$ is a basis of $Ker(\omega)$. One can easily check that this map is well defined.

\subsection{Affine coordinates on the family of $\alpha$-plane.}
Since the family of $\alpha$-planes is biholomorphic to $\CP^3$ we know that there is an affine covering $\{\osu_0,\ldots,\osu_3\}$ of the family of $\alpha$-planes. Let us write down the affine charts on $\osu_0$ and $\osu_1$. Let $v=(v_0,v_1,v_2,v_3)\in\C^4$ be a non-zero vector and let us assume that $v_0\ne0$, resp. $v_1\ne0$. Let $w_0:=v_0^{-1}v=(1,\zeta_1,\zeta_2,\zeta_3)$, resp. $w_1:=v_1^{-1}v=(\rho_1,1,\rho_2,\rho_3)$. Then the null 3-plane $\langle w_0\wedge f_1,w_0\wedge f_2,w_0\wedge f_3\rangle$, resp. $\langle w_1\wedge f_0,w_1\wedge f_2,w_1\wedge f_3\rangle$ has a unique basis of the form
\begin{equation}\label{1.affine coor on alpha planes}
\left(
\begin{matrix}
1&0&0\\
0&1&0\\
0&0&1\\
0&-\zeta_3&\zeta_2\\
\zeta_3&0&-\zeta_1\\
-\zeta_2&\zeta_1&0
\end{matrix}
\right), \ resp.\ 
\left(
\begin{array}{ccc}
-1&0&0\\
-\rho_2&\rho_1&0\\
-\rho_3&0&\rho_1\\
0&-\rho_3&\rho_2\\
0&0&-1\\
0&1&0
\end{array}
\right)
\end{equation}
where we use the notation from \ref{notation}. Similarly for $v_3\ne0$ and $v_4\ne0$. The change of coordinates between $\osu_0$ and $\osu_1$ is
\begin{eqnarray}\label{transition funct on cp3}
\zeta^{-1}_1=\rho_1,\zeta_2\zeta_1^{-1}=\rho_2,\zeta_3\zeta_1^{-1}=\rho_3.
\end{eqnarray}


\subsection{The set $\osw$.}  Let $x\in\osu$ and let $g\in\G_-$ be the unique element such that $\pi(g)=x$. Let $x^c:=g(x_0^c)$. Then $x^c$ is complement of $x$ in $x^\bot$ with basis $\{g(e_3),g(e_4),g(e_5),g(\bar e_3),g(\bar e_4),g(\bar e_5)\}$. Thus for any $x\in\osu$ we have the preferred isomorphisms $x^c\cong x_0^c\cong\C^6$. Any null 5-plane $z'\in\osw$ has a null orthogonal basis $\{v_1,\ldots,v_5\}$ such that $x':=\langle v_1,v_2\rangle\in\osu$ and that the null 3-plane $y':=\langle v_3,v_4,v_5\rangle$ belongs $(x')^c$. The null 3-plane $y'$ belongs to the family of $\alpha$-planes in $(x'_0)^c$. If $y'\in\osu_0$, resp. $y'\in\osu_1$ then $z'$ admits the unique basis of the form
\begin{equation}\label{coord on w0}
\left(
\begin{array}{ccccc}
1&0&0&0&0\\
0&1&0&0&0\\
0&0&1&0&0\\
0&0&0&1&0\\
0&0&0&0&1\\
z_{11}&z_{21}&0&-\zeta_3&\zeta_2\\
z_{21}&z_{22}&\zeta_3&0&-\zeta_1\\
z_{31}&z_{32}&-\zeta_2&\zeta_1&0\\
0&z_0&-z_{11}&-z_{21}&-z_{31}\\
-z_0&0&-z_{21}&-z_{22}&-z_{32}\\
\end{array}
\right),\ resp. \ \ \
\left(
\begin{array}{ccccc}
1&0&0&0&0\\
0&1&0&0&0\\
0&0&1&0&0\\
w_{31}&w_{32}&\rho_2&\rho_1&0\\
w_{21}&w_{22}&\rho_3&0&-\rho_1\\
w_{11}&w_{21}&0&-\rho_3&-\rho_2\\
0&0&0&0&1\\
0&0&0&1&0\\
0&w_0&\ast&\ast&\ast\\
-w_0&0&\ast&\ast&\ast\\
\end{array}
\right)
\end{equation}
where $\ast$ are determined by the other entries. Similarly for $y'\in\osu_2$ or $y'\in\osu_3$. Let us denote the corresponding affine subsets of $\osw$ by $\osw_i$ where $i=0,1,2,3$. Thus $\mathfrak W:=\{\osw_i|i=0,1,2,3\}$ is affine covering of $\osw$. We write the left hand side of (\ref{coord on w0}) in the block form 
\begin{equation}\label{block form on w0}
\left(
\begin{matrix}
1_2&0\\
0&1_3\\
B_1&B_2\\
B_0&-B_1^T\\
\end{matrix}
\right), 
B_0=\left(\begin{matrix}
0&z_0\\
-z_0&0
\end{matrix}
\right),B_1=(z_{ij})_{i=1,2,3}^{j=1,2},B_2=\left(\begin{matrix}
0&-\zeta_3&\zeta_2\\
\zeta_3&0&-\zeta_1\\
-\zeta_2&\zeta_1&0
\end{matrix}
\right).
\end{equation}
The change of coordinates on $\osw_0\cap\osw_1$ is 
\begin{eqnarray}\label{transition funct on w01}
w_0&=&z_0+z_{32}z_{21}\zeta_1^{-1}-z_{22}z_{31}\zeta_1^{-1},w_{11}=z_{11}+z_{31}\zeta_3\zeta_1^{-1}+z_{21}\zeta_2\zeta_1^{-1}\\
w_{12}&=&z_{12}+z_{32}\zeta_3\zeta_1^{-1}+z_{22}\zeta_2\zeta_1^{-1}\nonumber,w_{ij}=(-1)^i\zeta^{-1}_1 z_{ij},i=2,3,j=1,2\nonumber
\end{eqnarray}
and those in (\ref{transition funct on cp3}). 

\section{Sections of the bundle $\sO_\lambda$ over the set $\osw$.}\label{bdle on twistor space}
Let $\C_\lambda$ be an one-dimensional $\LGR$-module with highest weight $(\frac{5}{2},\frac{5}{2},\frac{5}{2},\frac{5}{2},\frac{5}{2})$. Let $\sO_\lambda$ be the sheaf of holomorphic sections of the bundle $\vbc:=\G\times_\LGR\C_\lambda$. Let $\osm\subset\G/\LGR$ be an open subset. Let $\sO_\lambda(\osm)$ be the space of holomorphic sections of the bundle $\vbc$ over $\osm$ and let $\sO(\osm)$ be the space of holomorphic functions on $\osm$. 

Let $z\in\osw_0$ and let $\{v_1,\ldots,v_5\}$ be the preferred basis of $z$ from (\ref{coord on w0}). Then $\{v_1,\ldots,v_5,\bar e_3,\bar e_4,\bar e_5,\bar e_1,\bar e_2\}$ is null orthonormal basis of $\C^{10}$. Let $g\in End_\C(\C^{10})$ be the linear map such that for all $1\le i\le 5:g(e_i)=v_i,g(\bar e_i)=\bar e_i$. Then $g\in\G$ and the map $z\mapsto\rho_0(z):=g$ is section of the principal $\LGR$-bundle over $\osw_0$. We define similarly for $i=1,2,3$ sections $\rho_i$ over $\osw_i$. Let us write transition function on $\vbc$ between $\osw_0$ and $\osw_1$ in the preferred trivializations $\rho_0$ and $\rho_1$. Let $f_0\in\sO(\osw_0),f_1\in\sO(\osw_1)$. Then $f_0,f_1$ defines an element of $\sO_\lambda(\osw_0\cup\osw_1)$ iff 
\begin{equation}\label{transition function for sections}
f_0(z_0,z_{ij},\zeta_i)=\zeta_1^{-5}f_1(w_0,w_{ij},\rho_i)
\end{equation}
on $\osw_0\cap\osw_1$. The transition functions on $\vbc$ between all sets $\osw_i$ are rational and thus we can also consider rational sections of $\vbc$ over $\osw$.

\subsection{Cohomology groups.} The Penrose transform gives an isomorphism 
\begin{equation}\label{penrose isomorphism}
\pt:H^3(\osw,\sO_\lambda)\cong Ker(D_1,\osu)
\end{equation}
between the sheaf cohomology of $\sO_\lambda$ over the set $\osw$ and the kernel of $D_1$ on $\osu$ consisting of holomorphic sections where the bundles from (\ref{sequence of operators}) are $\mathcal V_1=\G\times_\LGP(\C_\nu\otimes\Sp_+)$, $\mathcal V_2=\G\times_\LGP(\C_\mu\otimes\Sp_+)$ where $\GL(2,\C)$-module $\C_\nu$, resp. $\C^2_\mu$ has highest weight $(\frac{5}{2},\frac{5}{2})$, resp. $(\frac{7}{2},\frac{5}{2})$ and $\Sp_+\cong\C^4,\Sp_-\cong(\C^4)^\ast$ as $\SL(4,\C)$-modules. 

The Leray theorem states that $H^i(\osw,\sO_\lambda)\cong\check{H}^i(\mathfrak W,\sO_\lambda)$ where $\check{H}^\ast(\mathfrak W,\sO_\lambda)$ are the cohomology groups computed with respect to the affine covering $\mathfrak W$, see for example \cite{WW}. The co-chains groups are $C^4(\mathfrak W,\sO_\lambda)=0,C^3(\mathfrak W,\sO_\lambda)=\big\{(\bigcap_{i=0}^3\osw_i,f)|f\in\sO_\lambda(\bigcap_{i=0}^3\osw_i)\big\},C^2(\mathfrak W,\sO_\lambda)=\ldots$.
By definition we have that $\check H^3(\mathfrak W,\sO_\lambda):=C^3(\mathfrak W,\sO_\lambda)/Im(\delta^2)$ where $\delta^2$ is the \v{C}ech co-differential. We will denote the cohomology classes by $[\ ]$. Let us make some simple observations about $\check{H}^3(\mathfrak W,\sO_\lambda)$. We will work with the affine chart on $\osw_0$.

Let us first notice that $\bigcap_{i=0,1,2,3}\osw_i=\{z\in\osw_0|\zeta_1\ne0,\zeta_2\ne0,\zeta_3\ne0\}$. Thus if $f\in\sO_\lambda(\bigcap_{i=0,1,2,3}\osw_i)$ is a holomorphic section then $f$ is the converging sum of rational sections $f(s_0,s_{ij},r_k)=z_0^{s_0}\prod_{ij}z_{ij}^{s_{ij}}\zeta_1^{-r_1}\zeta_2^{-r_2}\zeta_3^{-r_3}$ where $s_0,s_{ij}\ge0$ and $r_1,r_2,r_3\in\Z$. It is easy to see that $[f(s_0,s_{ij},r_k)]=0$ if $r_1<0$ or $r_2<0$ or $r_3<0$. From the formula (\ref{transition function for sections}) follows that if $5+s_0+\sum s_{ij}>r_1+r_2+r_3$ then $f(s_0,s_{ij},r_k)$ extend to a rational section on $\osw_1\cap\osw_2\cap\osw_3$ and thus $[f(s_0,s_{ij},r_k)]=0$. However notice that this does not characterize $\check H^3(\mathfrak W,\sO_\lambda)$. For example the cohomology class of the section $z_{31}^i\zeta_1^{-1}\zeta_2^{-1}\zeta_3^{-3}$ is trivial for any $i\ge0$ although the relation does not hold. The full characterization of the 3-rd sheaf cohomology group of rational sections will be given in Theorem \ref{decomposition}. where we give it as a direct sum of $\G_0$-modules.

\section{The correspondence $x\in\osu\mapsto\eta\circ\tau^{-1}(x)\subset\osw$.} Let us write the correspondence on $\osw_0$. Then
\begin{eqnarray}\label{fiber on twistor space}
\left(
\begin{array}{c}
1_2\\
X_1\\
X_2\\
X_{12}-\frac{1}{2}(X_1^TX_2+X_2^TX_1)\\
\end{array}
\right)\in\osu\mapsto
\left(
\begin{array}{cc}
1_2&0\\
0&1_3\\
X_2-\zeta X_1&\zeta\\
X_{12}+\frac{1}{2}(X_2^TX_1-X_1^TX_2)+X_1^T\zeta X_1&-X^T_2-X_1^T\zeta\\
\end{array}
\right)\subset\osw_0
\end{eqnarray}
where the blocks of the matrix on the right hand side are those from (\ref{block form on w0}). In particular $X_1,X_2\in M(3,2,\C),X_{12}\in A(2,\C),\zeta\in A(3,\C)$. Notice that the fibres of the correspondence are biholomoprhic to $\CP^3$.

\subsection{Integration.}
The isomorphism (\ref{penrose isomorphism}) is given by integrating over the fibres of the correspondence (\ref{fiber on twistor space}). Let $f\in C^3(\osw,\sO_\lambda)$. We write $f\equiv f(B_0,B_1,B_2)$ where $(B_0,B_1,B_2)$ are the matrices from (\ref{block form on w0}). Then the integral formula in the affine chart on $\osw_0$ is 
\begin{eqnarray}\label{integrating over fibers}
\pt(f)_A(X_{12},X_1,X_2)=\frac{1}{(2\pi i)^3}\int_{(S^1)^3}(1,\zeta_1,\zeta_2,\zeta_3)f(X_{12}+\frac{1}{2}(X_2^TX_1-X_1^TX_2)+X_1^T\zeta X_2,X_2-\zeta X_1,\zeta)d\zeta_1d\zeta_2d\zeta_3
\end{eqnarray}
where $d\zeta_1d\zeta_2d\zeta_3$ is the holomorphic top form on the fibres which is homogeneous of degree $4$ in the homogeneous coordinates. The section $f$ is homogeneous of degree $-5$ in the homogeneous coordinates. Thus the integrand is homogeneous of degree zero and the integration does not depend on the choice of trivialization of the fibres of the correspondence. For example 
$$\pt(\zeta^{-1}_1\zeta_2^{-1}\zeta_3^{-1})=\frac{1}{(2\pi)^3}\int_{(S^1)^3}\frac{(1,\zeta_1,\zeta_2,\zeta_3)d\zeta_1d\zeta_2d\zeta_3}{\zeta_1\zeta_2\zeta_3}=(1,0,0,0)$$ is constant spinor on $\osu$. 

\section{Decomposition of monogenic sections into irreducible $\G_0$-modules.}
It is convenient to introduce a gradation on the space of all polynomial spinors on $\osu$. We trivialize the $\LGP$-bundle over $\osu$ by the map $\osu\xrightarrow{\pi^{-1}}\G_-\hookrightarrow\G$. We write coordinates on $\lag_-$ and thus also on $\osu$ as 
$X_1=(x^1_{ij})^{i=1,2,3}_{j=1,2},X_2=(x^2_{ij})^{i=1,2,3}_{j=1,2},
X_{12}=\left(\begin{matrix}
0&x_{12}\\
-x_{12}&0
\end{matrix}
\right).$ We will denote polynomials on the affine space $\lag_-$ with the same letters as the coordinates. 

Let us first define \textit{degree} of linear polynomials by setting $deg(x_{12}):=2,deg(x_{ij}^k):=1$. Let us extend it to the set of monomials in $\C[x_{12},x_{ij}^k]$ by requiring that $deg$ is morphism of $(\C[x_{12},x_{ij}^k],.)$ and $(\Z,+)$. Let $N_k$ be the vector space generated by the monomials of the degree $k$. Finally, let $M_k$ be the space of monogenic spinors whose components belong to $N_k$. If we extend this gradation naturally also to $\Gamma(\mathcal V_2)$ over $\osu$, then the operator $D_1$ is homogeneous of degree $-1$.

The space of linear monogenic spinors is $\G_0$-irreducible. The space of quadratic monogenic spinors $M_2$ decompose into
$\mW_{(\frac{9}{2},\frac{5}{2})}\otimes\mV_{(3,2,0,0)}\oplus\mW_{(\frac{7}{2},\frac{7}{2})}\otimes\mV_{(3,1,1,0)}\oplus\mW_{(\frac{7}{2},\frac{7}{2})}\otimes\mV_{(1,0,0,0)}$ where $\mW_{(a,b)}$, resp. $\mV_{(a,b,c,d)}$ stands for irreducible $\GL(2,\C)$, resp. $\SL(4,\C)$-module with highest weight $(a,b)$, resp. $(a,b,c,d)$. Let us write the dimensions of the modules, highest weight vectors from $\check H^3(\mathfrak W,\sO_\lambda)$ and the corresponding monogenic spinors. We have 
\begin{eqnarray*}
&&\mW_{(\frac{9}{2},\frac{5}{2})}\otimes\mV_{(3,2,0,0)}:180,\frac{z_{11}^2}{\zeta_1\zeta_2\zeta_3},((x^2_{11})^2,0,0,0)\\
&&\mW_{(\frac{7}{2},\frac{7}{2})}\otimes\mV_{(3,1,1,0)}:36,\frac{z_{11}z_{22}-z_{12}z_{21}}{\zeta_1\zeta_2\zeta_3},(x^2_{11}x^2_{22}-x^2_{21}x^2_{12},0,0,0)\\
&&\mW_{(\frac{7}{2},\frac{7}{2})}\otimes\mV_{(1,0,0,0)}:4,\frac{z_0}{\zeta_1\zeta_2\zeta_3}-\frac{z_{22}z_{31}-z_{21}z_{32}}{\zeta_1^2\zeta_2\zeta_3}-\frac{z_{11}z_{32}-z_{12}z_{31}}{\zeta_1\zeta_2^2\zeta_3}-\frac{z_{12}z_{21}-z_{11}z_{22}}{\zeta_1\zeta_2\zeta^2_3},\nonumber\\
&&(3x_{12}+\frac{1}{2}\sum_{i=1}^3(x^1_{i1}x^2_{i2}-x^2_{i1}x^1_{i2}),x^2_{21}x^2_{32}-x^2_{31}x^2_{22},x^2_{31}x^2_{12}-x^2_{11}x^2_{32},x^2_{11}x^2_{22}-x^2_{21}x^2_{12}).\nonumber
\end{eqnarray*}
In general we have the following theorem.

\begin{thm}\label{decomposition}
Let us keep the notation as above. Then the space $M_k$ of monogenic spinors of degree $k$ on $\osu$ decomposes into irreducible $\G_0$-modules 
\begin{equation}
M_k\cong\bigoplus_{a,b,l\ge0,2a+b+2l=k}\mW_{(\frac{5}{2}+l+a+b,\frac{5}{2}+l+a)}\otimes\mV_{(2a+b+1,a+b,a,0)}.
\end{equation}
The decomposition of algebraic monogenic spinors into irreducible $\G_0$-modules is multiplicity free.
\end{thm}

Proof: Let us recall that the action of $\G_0$ on $\check H^3(\mathfrak W,\sO_\lambda)$ is induced by the left action on the total space of the parabolic geometry. Let us now compute the weight of 
\begin{equation}\label{function1}
f=\frac{z_0^{s_0}\prod_{ij} z_{ij}^{s_{ij}}}{\zeta_1^{r_1}\zeta_2^{r_2}\zeta_3^{r_3}}\in C^3(\osw,\sO_\lambda).
\end{equation}
Let us denote for $j=1,2: c_j:=s_{1j}+s_{2j}+s_{3j},$ for $i=1,2,3:s_i:=s_{i1}+s_{i2}$ and let $r:=r_1+r_2+r_3,s:=s_1+s_2+s_3.$ We write weights as $\lagl(2,\C)\oplus\lasl(4,\C)$-weights. We find that the weight of $f$ is  
\begin{equation}\label{g0 weight}
(c_1+s_0+\frac{5}{2},c_2+s_0+\frac{5}{2})\oplus(5+s-r,r_1+s_1,r_2+s_2,r_3+s_3).
\end{equation}
Let us recall (\ref{isom sl(4,C) and so(6,C)}). Let us write the $\LGR$-module structure on $\C_\lambda$ by $\sigma$. Let $A_{12}$ be a standard positive root in $\lasl(2,\C)$. We find that all simple roots beside $E_{12}$ preserve the image of the section\footnote{See section \ref{bdle on twistor space}.} $\rho_0$ of the principal $\LGR$-bundle over $\osw_0$. Thus $f$ behaves as a rational function when we differentiate with simple roots other then $E_{12}$. For the root $E_{12}$ we have that $E_{12}f(z)=\frac{d}{dt}|_0\sigma(r^{-1}(z,t))f(z'(t))$ such that $z'(t)\in\osw_0$ and $\rho_0(z'(t))=\exp(-tE_{12})\rho_0(z)r(z,t)$ for a unique $r(z,t)\in\LGR$ for $t$ sufficiently small. If we differentiate we find out that $E_{12}f$ is the usual derivation of rational function which obeys the Leibniz rule plus the term $\dot\sigma(E_{12})f(z):=\frac{d}{dt}|_0\sigma(r^{-1}(t,z))f(z)$. We find that
\begin{eqnarray}\label{table1}
A_{12}z_{i2}&=&z_{i1},E_{34}\zeta_2=-\zeta_3,E_{34}z_{3i}=z_{2i},E_{23}\zeta_1=-\zeta_2, E_{23}z_{2i}=z_{1i},E_{43}\zeta_3=-\zeta_2, E_{43}z_{2i}=z_{3i},\nonumber\\
E_{32}\zeta_2&=&-\zeta_1,E_{32}z_{1i}=z_{2i},E_{21}\zeta_1=1,E_{12}z_{1i}=-\zeta_2z_{2i}-\zeta_3z_{3i},\nonumber\\
E_{12}z_0&=&z_{22}z_{31}-z_{21}z_{32},\dot\sigma(E_{12})f=5\zeta_1f, for\  i=2,3:E_{12}\zeta_i=\zeta_i\zeta_1,E_{12}z_{ij}=z_{ij}\zeta_1.
\end{eqnarray}
while all other terms are zero. This and Leibniz rule allows us to compute easily the action of $\lag_0$ on $C^3(\mathfrak W,\sO_\lambda)$.

\begin{lemma}\label{help3}
Let $f$ be a $\laq_0$-highest weight vector in the space of polynomials on the block $B_1$ in $(\ref{block form on w0})$, i.e. $f\in\C[z_{ij}]$. Then $f$ is $A(z_{11}z_{22}-z_{21}z_{12})^az_{11}^b$ for some $A\in\C,a,b=0,1,2,\ldots$
\end{lemma}
Proof: See \cite{GW}. $\Box$

\begin{lemma}\label{help1}
Let $f$ be the rational section from $(\ref{function1})$ such that the weight of $f$ is dominant and $s_0=0,r_i\ge 1$. Then the class  $[E_{12}^{r_1+r_2+r_3-3}E_{23}^{r_2+r_3-2}E_{34}^{r_3-1}f]\in \check H^3(\mathfrak W,\sO_\lambda)$ is non-zero.
\end{lemma}

Proof: We have that 
\begin{eqnarray*}
E_{34}^{r_3-1}f&=&A\frac{\prod_{ij}z_{ij}^{s_{ij}}}{\zeta_1^{r_1}\zeta_2^{r_2+r_3-1}\zeta_3}+\zeta_3^{-2}(\ldots)\\
E_{23}^{r_2+r_3-2}E_{34}^{r_3-1}f&=&B\frac{\prod_{ij}z_{ij}^{s_{ij}}}{\zeta_1^{r_1+r_2+r_3-2}\zeta_2\zeta_3}+\zeta^{-2}_2(\ldots)+\zeta_3^{-2}(\ldots)\\
E_{12}^{r_1+r_2+r_3-3}E_{23}^{r_2+r_3-2}E_{34}^{r_3-1}f&=&C\frac{\prod_{ij}z_{ij}^{s_{ij}}}{\zeta_1\zeta_2\zeta_3}+\zeta_1^{-2}(\ldots)+\zeta^{-2}_2(\ldots)+\zeta_3^{-2}(\ldots),\nonumber
\end{eqnarray*}
where $\ldots$ denotes sections where $\zeta_i$ appear only in denominators and
\begin{eqnarray*}
A&=&(-1)^{r_3-1}r_2(r_2+1)\ldots(r_2+r_3-2)\\
B&=&(-1)^{r_2+r_3-2}r_1(r_1+1)\ldots(r_1+r_2+r_3-3)A\\
C&=&AB(s_2+s_3+5-r)(s_2+s_3+5-(r-1))\ldots(s_2+s_3+1).
\end{eqnarray*}
Since the weight of $f$ is by assumption dominant, then $5+s-r\ge r_1+s_1>0$ and thus $5+s_2+s_3-r>0$. It follows that $C\ne0$ and thus also 
\begin{equation}\label{help4}
\pt(E_{12}^{r_1+r_2+r_3-3}E_{23}^{r_2+r_3-2}E_{34}^{r_3-1}f)=C(\prod_{ij}(x^2_{ij})^{s_{ij}},0,0,0)+\ldots
\end{equation}
where $\ldots$ denotes some spinors whose first components are different from $\prod_{ij}(x^2_{ij})^{s_{ij}}$. In particular we get that the cohomology class is non-zero.$\Box$

\begin{lemma}\label{help2}
Let 
\begin{equation}\label{function2}
f=\sum_{k=1}^K g_k,\ where\ g_k=\frac{f_k}{\prod_i\zeta_i^{r_i^k}},
\end{equation}
be a highest weight vector in $H^3(\osw,\sO_\lambda)$ such that all $f_k\in\C[z_{ij}]$. Then $K=r^1_1=r^1_2=r^1_3=1$ and $f_1$ is a $\laq_0$-maximal polynomial given in Lemma \ref{help3}.
\end{lemma}
Proof: Each summand in (\ref{function2}) satisfy the assumptions of Lemma (\ref{help1}). Let us notice that from (\ref{g0 weight}) follows that for all $1\le j,k\le K:r_1^j+r_2^j+r_3^j=r_1^k+r_2^k+r_3^k$ and that $deg(f_j)=deg(f_k)$. We can choose in (\ref{function2}) indexation by $j$ such that for all $k>1$ the following holds: $r_1^1>r_1^k$ or $r_1^1=r_1^k$ and $r_2^1>r_2^k$ or $r_1^1=r_1^k$ and $r_2^1=r_2^k$ and $r_3^1>r_3^k$. 

Let us assume that $r^1_1r^1_2r^3_3\ge2$. Let $E:=E_{12}^{r_1^1+r_2^1+r_3^1-3}E_{23}^{r_2^1+r_3^1-2}E_{34}^{r_3^1-1}.$ The formula (\ref{help4}) reveals that $\pt(E(g_1))\ne0.$ Similar manipulations give that $\pt(E(g_1))\ne -\pt(E(f-g_1)$ and thus $\pt(E.f)\ne0$ and thus $f$ is not highest weight vector. Thus the only possibility is that $K=r_1^1=r_2^1=r_3^1=1$ and $f_1$ is $\laq_0$-maximal. $\Box$

\begin{lemma}
Let 
\begin{equation}\label{general hwv}
f=\sum_{i=0}^{s_0}z_0^{s_0-i}f_i=z_0^{s_0}f_0+z_0^{s_0-1}f_1+\ldots
\end{equation} be a maximal highest weight vector such that $f_i\in C^3(\mathfrak W,\sO_\lambda)$ are rational sections which do not depend on $z_0$ and $[f_0]\ne0$. Then $f_0$ is also highest weight vector and $f$ is uniquely determined by $s_0$ and $f_0$. Conversely given a non-zero highest weight vector $f_0\in C^3(\mathfrak W,\sO_\lambda)$ that does not depend on the variable $z_0$ and $s_0\ge0$, then there exists a unique highest weight vector $f$ of the form as in (\ref{general hwv}) for some $f_i,i=1,\ldots,s_0$.
\end{lemma}
Proof: We easily check that if $f$ is highest weight vector then also $f_0$ is highest weight vector. Thus $f_0$ is a multiple non-zero of $z_{11}^a(z_{11}z_{22}-z_{12}z_{21})^b\zeta^{-1}_1\zeta^{-1}_2\zeta^{-1}_3$ for some $a,b\ge 0$. Let us check uniqueness of $f$ given $f_0$ and $s_0$. Let $\hat f,\tilde f$ be two highest weight vectors of the same weight such that $\hat f=z_0^{s_0}f_0+z_0^{s_0-1}(\ldots)$ and $\tilde f=z_0^{s_0}f_0+z_0^{s_0-1}(\ldots)$. Then $\check f:=\hat f-\tilde f=z_0^{t_0}\check f_0+z_0^{t_0-1}\check f_1+\ldots$ with $t_0<s_0$ has to be a highest weight vector with the same weight as $f$ and $f'$. Thus $\check f_0$ is a multiple of $z_{11}^c(z_{11}z_{22}-z_{12}z_{21})^d\zeta^{-1}_1\zeta^{-1}_2\zeta^{-1}_3$ for some $c,d\ge 0$. But the formula (\ref{g0 weight}) shows that then $a=c,b=d$ and thus also $t_0=s_0$. Contradiction.

Let us consider a filtration $\{F_i|i\ge0\}$ of $C^3(\mathfrak W,\sO_\lambda)$ given by the degree of $z_0$, i.e. $F_i:=\{g\in C^3(\mathfrak W,\sO_\lambda)|\partial_{12}^{i+1}g=0\}$ where $\partial_{12}$ is the coordinate vector field corresponding to the variable $z_0$. From the table (\ref{table1}) follows that $\G_0$ preserve this filtration. Let $f_{top}=z_0^{s_0}f_0$ be the highest part of $f$. Let $V=\G_0.f_{top}=\{\sum_{j=1}^Mg_j.f_{top}|g_j\in\G_0,M<\infty\}$. Then clearly $V$ is the smallest $\G_0$-module which contains the vector $f_{top}$. Moreover $V\subset F_{s_0}$ and from the table (\ref{table1}) follows that $V/F_{s_0-1}$ is spanned by $f_{top}$. We have that $V=\oplus_iV_i$ for some irreducible $\G_0$-modules $V_i$. Let $h_i$ be a maximal vector of $V_i$. Since the filtration $F_i$ is $\G_0$-equivariant, there exists $i$ such that $h_i=f_{top}+l.o.t.$ where $l.o.t.$ means lower order terms in $z_0$-variables. From the uniqueness we have that $h_i$ is up to a multiple the unique highest weight vector with the leading term $z_0^{s_0}f_0$. $\Box$

Thus we have that any highest weight vector is uniquely determined by its leading term $f_{top}$ with respect to the variable $z_0$. If $f_{top}=z_0^l(z_{11}z_{22}-z_{12}z_{21})^az_{11}^b$, then $f$ is highest weight vector of the module $W_{(\frac{5}{2}+a+b+s_0,\frac{5}{2}+s_0+b)}\otimes V_{(2a+b+1,a+b,a,0)}$. $\Box$

\end{document}